\documentclass[11pt]{article}
\usepackage{amsmath,amsthm,amsfonts,amssymb,mathrsfs,epsfig,bm}
\usepackage[usenames]{color}
\usepackage{cancel,xspace}

\parskip        \bigskipamount

\newtheorem{theorem}{Theorem}
\newtheorem{lemma}{Lemma}

\newtheorem{proposition}{Proposition}
\newtheorem{definition}{Definition}
\theoremstyle{remark}

\def\R{\mathbb{R}}
\def\P{\mathbb{P}}
\def\E{\mathbb{E}}

\def\FF{\mathscr{F}}

\renewcommand{\phi}{\varphi}
\renewcommand{\epsilon}{\varepsilon}
\newcommand{\1}{{\text{\Large $\mathfrak 1$}}}
\newcommand{\comp}{\raisebox{0.1ex}{\scriptsize $\circ$}}

\definecolor{mygray}{gray}{0.9}

\newcommand{\reg}{\operatorname{\mathcal{R}}}

\begin{document}

\title{\bf Integral representation of Skorokhod reflection}
\author{\sc Venkat Anantharam 
\thanks{Research supported by the ARO
MURI grant W911NF-08-1-0233, Tools for the Analysis and
Design of Complex Multi-Scale Networks, and by the NSF
grants CCF-0500234, CCF-0635372, CNS-0627161
and CNS-0910702.}  
\and  
\sc Takis Konstantopoulos
\thanks{Research supported by an EPSRC grant.}}
\date{\small \em 29 March 2010}
\maketitle

\begin{abstract}
We show that a certain integral representation of the 
one-sided Skorokhod reflection of a 
continuous bounded variation function
characterizes the reflection in that it 
possesses a unique maximal solution 
which solves the Skorokhod reflection problem.
\end{abstract}

\section{Introduction}
The Skorokhod reflection problem has a long history. Skorokhod \cite{SK61}
introduced it as a method for representing a diffusion process with
a reflecting boundary at zero.
Given a continuous function $X:[0,\infty) \to \R$, 
the standard Skorokhod reflection problem seeks to find 
$(Q(t), t \ge 0)$ and a 
continuous, nondecreasing function $Y: [0,\infty) \to \R_+$ with $Y(0)=0$,
such that $Q(t) := X(t) + Y(t) \ge 0$ for all $t$,
and $\int_0^\infty Q(s) dY(s)=0$.
Intuitively, the latter expresses the idea that $Y$ can increase
only at points $t$ such that $X(t)+Y(t)=0$.
Skorokhod \cite{SK61} showed that there is only one such $Y$, namely,
$Y(t) = -\inf_{0\le s \le t} (X(s) \wedge 0)$ and thus
\[
Q(t) = X(t) \vee \sup_{0 \le s \le t} (X(t)-X(s)).
\] 
We use the
standard notation $a \vee b := \max(a,b)$, $a \wedge b := \min(a,b)$.
The mapping $X \mapsto Q$ is referred to as the (one-sided) 
Skorokhod reflection mapping and has now become a standard 
tool in probability theory and other areas. As an example, 
we recall that if $X$ is the path of a Brownian motion
then $Q$ is a reflecting Brownian motion and $Q(t)$ has
the same distribution as $|X(t)|$ for all $t \ge 0$ \cite{CW,RY99}.
Several extensions of the Skorokhod reflection mapping exist
generalizing the range of $X$
(see, e.g., \cite{TAN79}) or its domain (see, e.g., \cite{AK05}).

The question resolved in this paper was motivated by an application
of the Skorokhod reflection in stochastic fluid queues \cite{KZD,KL}.
Suppose that $A, C$ are two jointly stationary and ergodic random
measures defined on a common probability space $(\Omega, \FF, \P)$,
with intensities $a, c$, respectively, such that $a < c$. Then there
exists a unique stationary and ergodic stochastic process
$(Q(t), t \in \R)$ defined on $(\Omega, \FF, \P)$
such that, for all $t_0 \in \R$, $(Q(t_0+t), t \ge 0)$ is
the Skorokhod reflection of $(Q(t_0)+A(t_0,t_0+t]-C(t_0,t_0+t],~ t \ge 0)$. 
In addition, if the random measures $A, C$ have no atoms then
\begin{equation}
\label{Qst}
Q(t) = \int_{-\infty}^t \1(Q(s) > C(s,t])~ dA(s),
\end{equation}
for all $t \in \R$, $\P$-almost surely.
The latter equation was called an ``integral representation'' of Skorokhod
reflection and extensions of it were formulated and proved in \cite{KL}.
The integral representation was found to be useful 
in several applications, e.g.\
(i) in deriving the so-called Little's law for stochastic fluid queues 
\cite{BB}, stating
that $\E [Q(0)] = (a/c) \E_A [Q(0)]$, where $\E_A$ is expectation with
respect to the Palm measure \cite{KALL} of $\P$ with respect to $A$,
and (ii) in deriving the form of the stationary distribution of 
a stochastic process derived from the local time of a L\'evy process
\cite{KKSS}.

In an open problems session of the workshop on
``New Topics at the Interface Between Probability and Communications''
\cite{K2010}, the second author asked whether and in what sense
\eqref{Qst} characterizes Skorokhod reflection.
The question will be made precise in Section \ref{theprob} below,
where the main theorem, Theorem \ref{mainthm},
which answers the question, is stated.
In Section \ref{directsec} the integral representation is 
explicitly proved, along with some auxiliary results.
Finally, in Section \ref{proofsec} a proof
of Theorem \ref{mainthm} is given.

\section{The problem}
\label{theprob}
Consider a locally finite signed measure $X$ on the Borel sets
of $\R$. Assume that $X$ has no atoms, i.e.\ $X(\{t\})=0$
for all $t \in \R$. 
Define
\begin{equation}
\label{Qstar}
Q^*(t) := \sup_{0 \le s \le t} X(s,t], \quad t \ge 0,
\end{equation}
where $X(s,t] = X((s,t])$ is the value of $X$ at the interval $(s,t]$.
\footnote{Since $X, A, C$ are assumed to have no atoms, we may as well write
$X[s,t]$ or $X(s,t)$ instead of $X(s,t]$, and likewise for $A$ and $C$,
 but we have chosen the notation to be
consistent with possible generalizations.}
In particular,
\[
Q^*(0)=0.
\]
Let $X(t) := X(0,t]$ and write \eqref{Qstar} as
\[
Q^*(t) = X(t) -\inf_{0 \le s \le t} X(s).
\]
The standard terminology \cite{CW,WHI} is
that $Q^*$ solves the Skorokhod reflection problem
for the function $t \mapsto X(t)$.

Decompose $X$ as the difference of two
locally finite nonnegative measures $A$, $C$,
without atoms, i.e.\ write
\begin{equation}
\label{Xdecomp}
X= A-C.
\end{equation}
We stress that $A$, $C$ are not necessarily the positive and negative parts
of $X$. In other words, the decomposition is not unique. For instance, 
we can add an arbitrary locally finite nonnegative measure without atoms
to both $A$ and $C$.

In \cite{KL} it was proved that \eqref{Qstar} also satisfies
the fixed point equation referred to as ``integral representation''
of the reflected process:
\begin{equation}
\label{ir}
Q(t) = \int_0^t \1(Q(s) > C(s,t])~ dA(s), \quad t \ge 0.
\end{equation}
A simpler version of this appeared earlier in \cite{KZD}; this
version was concerned with the case where $C$ is a multiple of
the Lebesgue measure.
In an open problems session of the workshop on
``New Topics at the Interface Between Probability and Communications''
\cite{K2010}, the second author asked whether and in what sense \eqref{ir} 
implies \eqref{Qstar}; the question was actually asked
for the special case where $C$ is a multiple of
the Lebesgue measure.

In this note we answer this question by proving the following:
\begin{theorem}
\label{mainthm}
Let $A$, $C$ be locally finite Borel measures on $\R_+ =[0,\infty)$
without atoms
and consider the integral equation \eqref{ir}.
This integral equation admits a unique maximal solution, i.e.\ a solution
which pointwise dominates any other solution. Further, this
maximal solution is precisely the function $Q^*$ defined by \eqref{Qstar}.
\end{theorem}

We proceed as follows. First, we present some auxiliary results
and also give a proof of \eqref{Qstar} $\Rightarrow$ \eqref{ir} which
is different from the one found in \cite{KL}. Then we prove Theorem
\ref{mainthm} by a successive approximation scheme
and by proving a number of lemmas.

\section{Proof of the integral representation and auxiliary results}
\label{directsec}
We first exhibit some properties of $Q^*$, defined by \eqref{Qstar},
and also show that
$Q^*$ satisfies the integral equation \eqref{ir}.
The proof of the latter in the special case where $C$ is a multiple
of the Lebesgue measure can be found in 
\cite[Lemma 1]{KZD} and in \cite[\S 3.5.3]{BB}.
A more general case is dealt with in \cite[Theorem 1]{KL}.
We give a different proof in Proposition \ref{direct} below.
The lemmas below are straightforward and well-known but we give proofs
for completeness.
As before, $X$ is a locally finite Borel measure without atoms
and $X=A-C$ is a
decomposition as the difference of two nonnegative locally finite
Borel measures without atoms.
We set
\[
A(t):= A(0,t], \quad
C(t):= C(0,t].
\]
\begin{lemma}
\label{itwontempty}
If $0 \le s \le s' \le t$ and if $Q^*(s) > C(s,t]$
then $Q^*(s') > C(s',t]$.
\end{lemma}
\proof
Assume that $C(s,t] < Q^*(s) = \sup_{0 \le u \le s} X(u,s]$.
This is equivalent to
\begin{align*}
C(t)-C(s) &< \sup_{0 \le u \le s} \{ A(s)-A(u) - (C(s)-C(u)) \}
\\
&= A(s) + \sup_{0 \le u \le s} \{-A(u)+C(u)\} -C(s),
\\
\text{that is, } \quad 
C(t) &< A(s) + \sup_{0 \le u \le s} \{-A(u)+C(u)\}.
\end{align*}
The right-hand side of the latter
is increasing in $s$ and so replacing $s$ by a larger $s'$ we
obtain
\[
C(t) < A(s') + \sup_{0 \le u \le s'} \{-A(u)+cu\} ,
\]
which is equivalent to $Q^*(s') > C(s',t]$.
\qed

\begin{lemma}
\label{semigroup}
$Q^*$ satisfies
\begin{equation}
\label{semigroupequ}
Q^*(t) = \sup_{s \le u \le t}
X(u,t] \vee (Q^*(s) + X(s,t]),
\quad 0 \le s \le t.
\end{equation}
\end{lemma}
\proof
We show that the right-hand side of \eqref{semigroupequ} equals 
the left-hand side.
\begin{align*}
\sup_{s \le u \le t} X(u,t] \vee (Q^*(s) + X(s,t] )
&= \sup_{s \le u \le t} X(u,t] \vee \{(\sup_{0 \le u \le s} X(u,s]) + X(s,t] \}
\\
&= \sup_{s \le u \le t} X(u,t] \vee  \sup_{0 \le u \le s}\{X(u,s]+X(s,t]\}
\\
&= \sup_{s \le u \le t} X(u,t]  \vee \sup_{0 \le u \le s} X(u,t]
\\
&= \sup_{0 \le u \le t} X(u,t] = Q^*(t).
\end{align*}

\begin{lemma}
\label{asinmypaper}
If $0 \le s \le t$ and  if $Q^*(s) \ge C(s,t]$ then
$Q^*(t) = Q^*(s) + X(s,t]$.
\end{lemma}
\proof
We use equation \eqref{semigroupequ}, rewritten as follows:
\begin{equation}
\label{QwithX}
Q^*(t) = \sup_{s \le u \le t} \big\{ X(u,t] \vee (Q^*(s) + X(s,t]) \big\}.
\end{equation}
Suppose $0\le s \le u \le t$ and that $Q^*(s) \ge C(s,t]$.
Then $Q^*(s) \ge C(s,u]$ and so
\begin{align*}
Q^*(s) + X(s,t] & \ge C(s,u] + X(s,t]
\\
&= C(s,u] + A(s,t]-C(s,t]
\\
& = A(s,t] - C(u,t]
\\
&\ge A(u,t] - C(u,t] = X(u,t],
\end{align*}
and this inequality implies that the term $X(u,t]$ inside the
bracket of the right-hand side of \eqref{QwithX} is not needed. Hence 
$Q^*(t) = Q^*(s) + X(s,t]$, which is what we wanted to prove.
\qed

Define next
\[
\sigma^*(t) := \sup\{0 \le s \le t:~ Q^*(s) \le C(s,t]\}.
\]
By Lemma \ref{itwontempty},
\begin{subequations}            
\label{beforeafter}
\begin{align}
Q^*(s) & \le C(s,t], \quad \text{ if } 0 \le s \le \sigma^*(t),
\label{before}\\
Q^*(s) & > C(s,t], \quad \text{ if } \sigma^*(t) < s \le t,
\label{after}
\end{align}
\end{subequations}
provided that the last inequality is non-vacuous.
Since the function $Q^*$ is nonnegative and continuous, we also have
\[
Q^*(\sigma^*(t)) = C(\sigma^*(t), t].
\]
\begin{proposition}
\label{direct}
If $X$ is a locally finite signed Borel measure on $[0,\infty)$ without atoms
and if $X=A-C$ is any decomposition of $X$ as the difference of
two nonnegative locally finite Borel measures without atoms, 
then the function $Q^*$
defined by \eqref{Qstar} satisfies \eqref{ir}.
\end{proposition}
\proof
By Lemma \ref{asinmypaper}, and the last display,
\begin{align*}
Q^*(t) &= Q^*(\sigma^*(t)) + A(\sigma^*(t), t] - C(\sigma^*(t), t]
\\
&= A(\sigma^*(t), t]
\\
&= \int_{\sigma^*(t)}^t dA(s)
\\
&= \int_0^t \1(Q^*(s) > C(s,t])~ dA(s),
\end{align*}
which is the integral representation formula \eqref{ir}.
Note that, to obtain the last equality in the last display,
we used \eqref{before}-\eqref{after}.
\qed

\section{Proof of Theorem \ref{mainthm}}
\label{proofsec}
{\em A priori}, it is not clear
that \eqref{ir} admits a maximal solution and, even if it does,
whether it satisfies \eqref{Qstar}.
We shall show the validity of these claims in the sequel.

We fix two locally finite measures $A$ and $C$ and 
define the map $\Theta$ on the set of nonnegative measurable functions by
\begin{equation}
\label{thetadef}
\Theta(Q)(t) :=  \int_0^t \1(Q(s) > C(s,t])~ dA(s), \quad t \ge 0.
\end{equation}
The integral equation \eqref{ir} then reads
\[
Q = \Theta(Q).
\]
We observe that $\Theta$ is increasing:
\begin{equation}
\label{IProp}
\text{ If $Q \le \widetilde Q$ then $\Theta(Q) \le \Theta(\widetilde Q)$.}
\end{equation}
Here, and in the sequel, given two functions $f, g :[0,\infty) \to \R$,
we write $f \le g$ to mean that $f(t) \le g(t)$
for all $t \ge 0$.
To see that \eqref{thetadef} holds,
simply observe that $Q \le \widetilde Q$ implies
$\1(Q(s) > C(s,t]) \le \1(\widetilde Q(s) > C(s,t])$ for
all $0\le s \le t$.

Define next a sequence of functions $(Q_k, k=0,1,2,\ldots)$ by
first letting
\[
Q_0 := \infty,
\]
and then, recursively,
\[
Q_{k+1} := \Theta(Q_k), \quad k \ge 0.
\]
Clearly, $Q_1(t) = \int_0^t dA(s) = A(t)$.
So $Q_0 \ge Q_1$. 
Since $\Theta$ is an increasing map, we see that, 
\[
Q_k \ge Q_{k+1} \ge 0 , \quad k \ge 0.
\]
We can then define
\[
Q_{\infty}(t) := \lim_{k \to \infty} Q_k(t).
\]
\begin{lemma}
\label{Qstarsmaller}
If $Q = \Theta(Q)$ then $Q \le Q_\infty$.
Furthermore,
\begin{equation*}
Q^* \le Q_\infty.
\end{equation*}
\end{lemma}
\proof
Suppose that $Q$ satisfies $Q=\Theta(Q)$.
Since the integrand in the right-hand side of \eqref{thetadef}
is $\le 1$, 
we have $Q(t) \le A(t)$ for all $t \ge 0$.
Letting $\Theta^{(k)}$ be the $k$-fold composition of $\Theta$ with itself, we
have
\[
Q = \Theta^{(k)}(Q) \le \Theta^{(k)}(A) = Q_k,
\]
and so $Q \le Q_\infty$.
In particular, Proposition \ref{direct} states that $Q^*=\Theta(Q^*)$.
Hence $Q^* \le Q_\infty$. 
\qed

However, it is not yet clear at this point that $Q_\infty$ is
a fixed point of $\Theta$.
We can only show that
\[
Q_\infty \ge \Theta(Q_\infty).
\]
Indeed,
$Q_\infty \le Q_k$ for all $k$, and so 
$\1(Q_\infty(s) > C(s,t]) \le \1(Q_k(s) > C(s,t])$, for all $0 \le s \le t$,
implying that $\Theta(Q_\infty) \le \Theta(Q_k) = Q_{k+1}$, and, by 
taking limits, that $\Theta(Q_\infty) \le Q_\infty$.

\begin{definition}[Regulating functions]
Consider functions $B:[0,\infty) \to [0,\infty)$ which are
continuous, nondecreasing, with $B(0)=0$, such that
$X(0,t] +  B(t) \ge 0$ for all $t \ge 0$.
Call these functions {\em regulating functions of $X$}.
The set of regulating functions is denoted by $\reg(X)$.
\end{definition}
We define a mapping 
\begin{align}
\label{Phi}
\Phi: & \reg(X) \to \reg(X)
\end{align}
in two steps:
Given $B \in \reg(X)$, first define 
\[
\sigma_B(t) := \sup\{0 \le s \le t:~ A(s)+B(s)-C(t) \le 0\},
\quad t \ge 0.
\]
Then let
\[
\Phi(B)(t) := B(\sigma_B(t)) ,
\quad t \ge 0.
\]
We actually need to show that what is claimed in \eqref{Phi} holds.
Namely: 
\begin{lemma}
If $B \in \reg(X)$ then $\Phi(B) \in \reg(X)$.
\end{lemma}
\proof
Clearly, $\sigma_B(\cdot)$ is nondecreasing. Since $B$ is nondecreasing,
it follows that $\Phi(B) = B \comp \sigma_B$ is nondecreasing.
Also, $\Phi(B)(0) = B(\sigma_B(0))=B(0)=0$. 
{From} the continuity of $A$, $B$ and the definition of $\sigma_B$, we
have
\begin{equation}
\label{sb}
A(\sigma_B(t))+B(\sigma_B(t))= C(t), \quad t \ge 0.
\end{equation}
We also have,
\begin{align*}
A(t) + \Phi(B)(t)-C(t) &= A(t) + B(\sigma_B(t))-C(t) \\
&= [A(t)-A(\sigma_B(t))] + [A(\sigma_B(t))+B(\sigma_B(t))-C(t)] \\
&= A(t)-A(\sigma_B(t)) \ge 0,
\end{align*}
where we used \eqref{sb} in the third step.
It remains to show that $\Phi(B)(\cdot)$ is continuous. 
Note that $\sigma_B(\cdot)$ need not be continuous.
However, $C(\cdot)$ is a continuous function and 
so, by \eqref{sb}, 
$t \mapsto A(\sigma_B(t)) + B(\sigma_B(t))$ is continuous. Hence 
\[
[A(\sigma_B(t+))-A(\sigma_B(t-)] + [B(\sigma_B(t+))-B(\sigma_B(t-))]=0,
\quad\text{for all $t$}.
\]
Since $A(\sigma_B(\cdot))$ and $B(\sigma_B(\cdot))$ are both nondecreasing, it
follows that $A(\sigma_B(t+))-A(\sigma_B(t-) \ge 0$ and $B(\sigma_B(t+))-B(\sigma_B(t-)) \ge 0$
and, since their sum is zero, they are both zero, implying
that $A(\sigma_B(\cdot))$ and $B(\sigma_B(\cdot))$ are continuous.
\qed

An immediate property of $\Phi$ is that
\begin{equation}
\label{PhiProp}
\Phi(B) \le B \quad \text{for all $B \in \reg(X)$}.
\end{equation}
Indeed, for all $t \ge 0$, $\sigma_B(t) \le t$ and so $B(\sigma_B(t)) \le B(t)$.

Starting with the function
\begin{equation}
\label{B1}
B_1(t):= C(t), \quad t \ge 0,
\end{equation}
we recursively define
\begin{equation}
\label{Bk}
B_{k+1} := \Phi(B_k), \quad k \ge 1.
\end{equation}
Therefore
\begin{equation}
\label{Binfty}
B_1 \ge B_2 \ge \cdots \ge B_k \downarrow B_\infty, \quad \text{
as } k \to \infty,
\end{equation}
where the inequalities and the limit are pointwise.
\begin{lemma}
\label{Binfreg}
The function $B_\infty$, defined via \eqref{B1}, \eqref{Bk} and \eqref{Binfty},
is a member of the class $\reg(X)$.
\end{lemma}
\proof
$B_\infty$ is nondecreasing since all the $B_k$ are nondecreasing.
Also, $B_\infty(0)=0$.
Since for all $k$, $A+B_k-C \ge 0$, we have
$A+B_\infty -C \ge 0$.
We proceed to show that $B_\infty$ is a continuous function.
We observe that, for $0 \le t \le t'$,
\begin{align*}
| \Phi(B)(t')-\Phi(B)(t)| 
&= |B(\sigma_B(t'))-B(\sigma_B(t))|
\\
&=  B(\sigma_B(t'))-B(\sigma_B(t))
\\
&\le A(\sigma_B(t'))-A(\sigma_B(t)) + B(\sigma_B(t'))-B(\sigma_B(t))
\\
&= [A(\sigma_B(t')) + B(\sigma_B(t'))] - [A(\sigma_B(t))+B(\sigma_B(t))]
\\
&= C(t')-C(t),
\end{align*}
where we again used \eqref{sb}.
It follows that the family of functions $\{\Phi(B), B \in \reg(X)\}$
is uniformly bounded and equicontinuous on each compact interval
of the real line. By the Arzel\`a-Ascoli theorem, the family
is compact and therefore $B_\infty$ is continuous.
We have established that $B_\infty \in \reg(X)$.
\qed

We now claim that $B_\infty$ is a fixed point of $\Phi$.
\begin{lemma}
\label{Binftyisfixed}
$\Phi(B_\infty)=B_\infty$.
\end{lemma}
\proof
By definition,
\[
\Phi(B_\infty)(t)= B_\infty(\sigma_{B_\infty}(t)),
\]
where
\[
\sigma_{B_\infty}(t) = \sup\{0 \le s \le t:~ A(s)+B_\infty(s) \le C(t)\}.
\]
Now, since $B_k \ge B_{k+1}$ for all $k \ge 1$, it follows
that $\sigma_{B_k} \le \sigma_{B_{k+1}}$ for all $k \ge 1$, and so
\[
\sigma_L(t) := \lim_{k \to \infty} \sigma_{B_k}(t)
\]
is well-defined.
Since $B_k \ge B_\infty$ for all $k \ge 1$, we have $\sigma_{B_k} \le
\sigma_{B_\infty}$.
Taking limits, we find 
\[
\sigma_L \le \sigma_{B_\infty}.
\]
Using the last two displays and the fact
that $B_k$ and $B_\infty$ are nondecreasing, we have
\begin{align*}
\Phi(B_\infty)(t) = B_\infty(\sigma_{B_\infty}(t))
& \ge B_\infty(\sigma_L(t))
\\
&= \lim_{k \to \infty} B_k(\sigma_L(t))
\\
&\ge \lim_{k\to \infty} B_k(\sigma_{B_k}(t))
\\
&= \lim_{k \to \infty} B_{k+1}(t) = B_\infty(t).
\end{align*}
By inequality \eqref{PhiProp}, $\Phi(B) \le B$ for all $B \in \reg(X)$ 
and since, by Lemma \ref{Binfreg},
$B_\infty \in \reg(X)$, it follows that we also have 
$B_\infty \le \Phi(B_\infty)$.
Therefore $B_\infty = \Phi(B_\infty)$, as claimed.
\qed

\begin{lemma}
Consider the function $Q^*$ defined by \eqref{Qstar} and
define a function $U$ by
\[
U(t) := Q^*(t)-X(0,t], \quad t \ge 0.
\]
Then
\begin{itemize}
\item[(i)]
$U \in \reg(X)$.
\item[(ii)]
$U = \Phi(U)$.
\end{itemize}
\end{lemma}
\proof
(i) We have $X(0,t]+U(t) = Q^*(t) \ge 0$ for all $t$.
Using \eqref{Qstar} and \eqref{Xdecomp} we see that
\begin{equation}
\label{Urep}
U(t) = \sup_{0 \le s \le t} \{-A(s)+C(s)\}.
\end{equation}
Therefore, $U(0)=0$, and $U$ is a continuous and nondecreasing.
We conclude that $U \in \reg(X)$.
To prove (ii), recall that $\Phi(U) = U \comp \sigma_U$ where
\[
\sigma_U(t) = \sup\{0 \le s \le t:~ A(s)+U(s) \le C(t)\}.
\]
Splitting the supremum in \eqref{Urep} in two parts, we obtain 
\begin{align*}
U(t) &
= \sup_{0\le s \le \sigma_U(t)} \{-A(s)+C(s)\}
\vee \sup_{\sigma_U(t) \le s \le t} \{-A(s)+C(s)\}.
\\
&= U(\sigma_U(t)) \vee \sup_{\sigma_U(t) \le s \le t} \{-A(s)+C(s)\}.
\end{align*}
For $s \ge \sigma_U(t)$, we have $A(s)+U(s) \ge C(t)$, i.e.\ 
$-A(s)+C(s) \le U(s)-C(s,t]$. Therefore
\begin{align*}
U(t) &\le U(\sigma_U(t)) \vee \sup_{\sigma_U(t) \le s \le t} \{U(s)-C(s,t]\}
\\
&= U(\sigma_U(t) = \Phi(U)(t).
\end{align*}
Thus, $U \le \Phi(U)$.
On the other hand, since $U \in \reg(X)$, we have $\Phi(U) \le U$,
by \eqref{PhiProp}.
\qed

\begin{lemma}
\label{BU}
Let $B \in \reg(X)$ be any fixed point of $\Phi$. Then $B \le U$.
\end{lemma}
\proof
Since $B=\Phi(B)=B\comp \sigma_B$ we have
\[
B = B \comp \sigma_B^{(k)}
\]
where
$\sigma_B^{(k)}:= \underbrace{\sigma_B \comp \cdots \comp \sigma_B}_{k \text{ times}}$.
Since
\[
t \ge \sigma_B(t) \ge \sigma_B \comp \sigma_B(t) \ge \cdots 
\ge \sigma_B^{(k)}(t),
\]
we may define
\[
\sigma_B^{(\infty)}(t) := \lim_{k \to \infty} \sigma_B^{(k)}(t).
\]
By the continuity of $B$,
\begin{equation}
\label{BBB}
B= B \comp \sigma_B^{(\infty)}.
\end{equation}
On the other hand, \eqref{sb} gives
\[
A \comp \sigma_B^{(k+1)} + B \comp \sigma_B^{(k+1)} = C \comp \sigma_B^{(k)}, \quad k \ge 1.
\]
Taking the limit as $k \to \infty$, and using the continuity
of $A$, $B$ and $C$, we have
\[
A \comp \sigma_B^{(\infty)} + B \comp \sigma_B^{(\infty)} = C \comp \sigma_B^{(\infty)}. 
\]
Since $A(t) + U(t) \ge C(t)$ for all $t$, we have
\[
A \comp \sigma_B^{(\infty)} + U \comp \sigma_B^{(\infty)} \ge C \comp \sigma_B^{(\infty)},
\]
and from the last two displays we conclude that
\[
U \comp \sigma_B^{(\infty)} \ge B \comp \sigma_B^{(\infty)}.
\]
Since $U$ is nondecreasing and since \eqref{BBB} holds, we have
\[
U \ge U \comp \sigma_B^{(\infty)} \ge B \comp \sigma_B^{(\infty)} = B,
\]
as claimed.
\qed

We are now ready to prove Theorem \ref{mainthm}.
We already know from Lemma \ref{Qstarsmaller} that $Q^* \le Q^\infty$.
So we only have to prove the opposite inequality.
Recall that $Q_1=A$ and $B_1 = C$.
Trivially then
\[
Q_1(t)+C(t) = A(t)+B_1(t), \quad t \ge 0.
\]
Thus, for $0 \le s \le t$ we have
\begin{align*}
Q_1(s) > C(s,t] 
& \iff
Q_1(s)+C(s) > C(t)
\\
& \iff
A(s)+B_1(s) > C(t)
\\
& \iff
s > \sigma_{B_1}(t).
\end{align*}
{From} this we get
\begin{align*}
Q_2(t) &= \int_0^t \1(Q_1(s) > C(s,t])~ dA(s) 
\\
&= \int_0^t \1(s > \sigma_{B_1}(t))~ dA(s)
\\
&= A(t)-A(\sigma_{B_1}(t)).
\end{align*}
But \eqref{sb} gives
\[
A(\sigma_{B_1}(t)) + B_1(\sigma_{B_1}(t)) = C(t),
\]
and so
\[
Q_2(t) + C(t) = A(t) + B_1(\sigma_{B_1}(t)) = A(t) + B_2(t), \quad t \ge 0.
\]
We now claim that
\[
Q_k(t) + C(t) = A(t) + B_k(t), \quad t \ge 0, \quad k \ge 1.
\]
This can be proved by induction along the same lines as above.
Taking limits as $k \to \infty$, we conclude
\[
Q_\infty(t) + C(t) = A(t) + B_\infty(t), \quad t \ge 0.
\]
Lemma \ref{Binftyisfixed} tells us that
$B_\infty$ is a fixed point of $\Phi$, and so, by
Lemma \ref{BU},
\[
B_\infty \le U.
\]
Hence
\begin{align*}
Q_\infty(t) + C(t) &= A(t) + B_\infty(t)
\\
&\le A(t)+U(t)
\\
&= Q^*(t)+C(t), \quad t \ge 0,
\end{align*}
and this gives 
\[
Q_\infty \le Q^*,
\]
as needed.
\qed

\section*{Acknowledgments}
We thank the Isaac Newton Institute for Mathematical Sciences for 
providing the stimulating research atmosphere
where this research work was done.

\vspace*{1cm}
\noindent
\hfill
\begin{minipage}[t]{6cm}
\small \sc
Authors' addresses:\\[3mm]
Venkat Anantharam\\
EECS Department\\
University of California\\
Berkeley, CA 94720, USA\\
E-mail: {\tt ananth@eecs.berkeley.edu}
\\[5mm]
Takis Konstantopoulos\\
School of Mathematical Sciences\\
Heriot-Watt University\\
Edinburgh EH14 4AS, UK\\
E-mail: {\tt takis@ma.hw.ac.uk}
\end{minipage}


{\small
\begin{thebibliography}{20}

\bibitem{AK05}
{\sc Anantharam, V. and Konstantopoulos, T.}
Regulating functions on partially ordered sets.
{\em Order} {\bf 22}, 145-183, 2005. 

\bibitem{BB}
{\sc Baccelli, F. and Br\'emaud, P.}
{\em Elements of Queueing Theory.}
Springer-Verlag, 2003.

\bibitem{CW}
{\sc Williams R. and Chung, K.-L.}
{\em An Introduction to Stochastic Integration.}
Birkh\"auser, Boston, 1989.

\bibitem{KALL}
{\sc Kallenberg, O.}
{\em Foundations of Modern Probability,} 2nd ed.
Springer-Verlag, New York, 2002.

\bibitem{KKSS}
{\sc Konstantopoulos, T., Kyprianou, A., Sirvi\"o, M., and Salminen, P.}
Analysis of stochastic fluid queues driven by local time processes.
{\em Adv.\ Appl.\ Probability} {\bf  40}, 1072-1103, 2008. 

\bibitem{KL}
{\sc Konstantopoulos, T., and Last, G.}
On the dynamics and performance of stochastic fluid systems.
{\em J.\ Appl.\ Prob.} {\bf 37}, 652-667, 2000. 

\bibitem{KZD}
{\sc Konstantopoulos, T., Zazanis, M. and de Veciana, G.}
Conservation laws and reflection mappings with an application 
to multiclass mean value analysis for stochastic fluid queues. 
{\em Stoch.\ Proc.\ Appl.} {\bf  65}, No. 1, 139-146, 1997.

\bibitem{K2010}
{\sc Konstantopoulos, T.}
Open problems session of the workshop on
{\em New Topics at the Interface Between Probability and Communications.}
Thursday, 14 January, 2010.

\bibitem{RY99}
{\sc Revuz, D.\ and Yor, M.}
{\em Continuous martingales and Brownian motion.}
Springer-Verlag, New York, 1999.

\bibitem{SK61}
{\sc Skorokhod, A.V.}
Stochastic equations for diffusions in a bounded region, 
{\em Theory Probab.\ Appl.} {\bf 6}, 264-274, 1961.

\bibitem{TAN79}
{\sc Tanaka, H.}
Stochastic differential equations with reflecting 
boundary condition in convex regions.
{\em Hiroshima Math.\ J.} {\bf 9}, 163-177, 1979.

\bibitem{WHI}
{\sc Whitt, W.}
{\em Stochastic-Process Limits.}
Springer-Verlag, New York, 2002.
\end{thebibliography}

}
\end{document}